\documentclass[10pt,a4paper,final,onecolumn]{IEEEtran}
\usepackage[utf8]{inputenc}
\usepackage{amsmath}
\usepackage{amsfonts}
\usepackage{amssymb}
\usepackage{graphicx}
\usepackage{epstopdf}
\usepackage{caption}
\usepackage[noadjust]{cite}
\usepackage{subcaption}
\usepackage{color}

\title{Method of Moments for computing electromagnetic scattering from a rough cylinder}
\author{Rahul Trivedi and~Uday K Khankhoje%
\thanks{R.~Trivedi and U.~K.~Khankhoje are with the Department of Electrical Engineering, Indian Institute of Technology Delhi, Hauz Khas, New Delhi 110016, India. E-mails: rtrivedi1995@gmail.com and uday@alumni.caltech.edu, respectively.}}%

\begin{document}
\maketitle
\begin{abstract}
In this tutorial paper, we formulate a two-dimensional integral-equation based method of moments approach for numerically computing the electromagnetic fields scattered from an azimuthally-rough dielectric cylinder or an axially-rough perfectly conducting cylinder. The electric field integral equation is discretized in both the cases to yield a system of linear equations, and a detailed evaluation of the matrix elements involved in this system of equations is presented. Numerical computation of important quantities such as the scattering cross-section or the scattered electric fields is also discussed. The primary purpose of this numerical scheme is to validate perturbative solutions that arise in the context of rough cylinder scattering problems.
\end{abstract}

\section{Introduction}

Evaluating the electromagnetic scattering from cylinders with rough surfaces is a computationally tedious task. This is because the stochastic nature of the surface necessitates taking an ensemble average of the scattered fields over many realizations of a rough surface till suitable numerical convergence is achieved. In spite of their accuracy, full-wave techniques such as integral-equation \cite{p} or finite-element \cite{fem} based methods, are impractical in remote-sensing applications where the scattering from a variety of cylinder configurations needs to be estimated \cite{smap,mariko}. Thus, it is natural to formulate perturbative solutions, which, for some tradeoff in accuracy, provide fast estimates of the scattering solution. 

A second order perturbative method for computing the electromagnetic fields scattered by a homogeneous cylinder with a rough surface has recently been proposed \cite{tap,piers}. In order to validate such a solution, it is essential to compare the perturbative solution with the solution predicted by a full-wave technique. In this paper, we describe the formulation of an integral-equation based method of moments which is used for the purpose of validating the above referenced perturbative solution. 

While the most general formulation of this problem is in three-dimensions, we restrict ourselves to two possible two-dimensional cylinder configurations -- the extension to three dimensions is straightforward, but beyond the scope of the present paper. In the first configuration, we consider scattering from a dielectric cylinder of infinite length whose roughness function has only azimuthal dependence. By employing the two-dimensional Green's function in a homogeneous medium, followed by an expansion of unknown currents in triangular basis functions and testing along pulse functions, we formulate our scattering solution in Section \ref{SecAzR}. In the second configuration, treated in Section \ref{SecAxR}, we consider scattering from a perfectly conducting cylinder of finite length with a roughness function which is only axially-dependent. Unlike the previous case, we require the three-dimensional Green's function here, and formulate the solution by expanding the unknowns in pulse basis functions and test along delta functions. Detailed remarks on the discretization parameters and evaluation of the associated matrix elements is made in both cases.

\section{Scattering from an azimuthally rough dielectric cylinder}\label{SecAzR}

\begin{figure}[htpb]
\centering
\includegraphics[trim = 30cm 13.7cm 33cm 2.8cm, scale = 0.77]{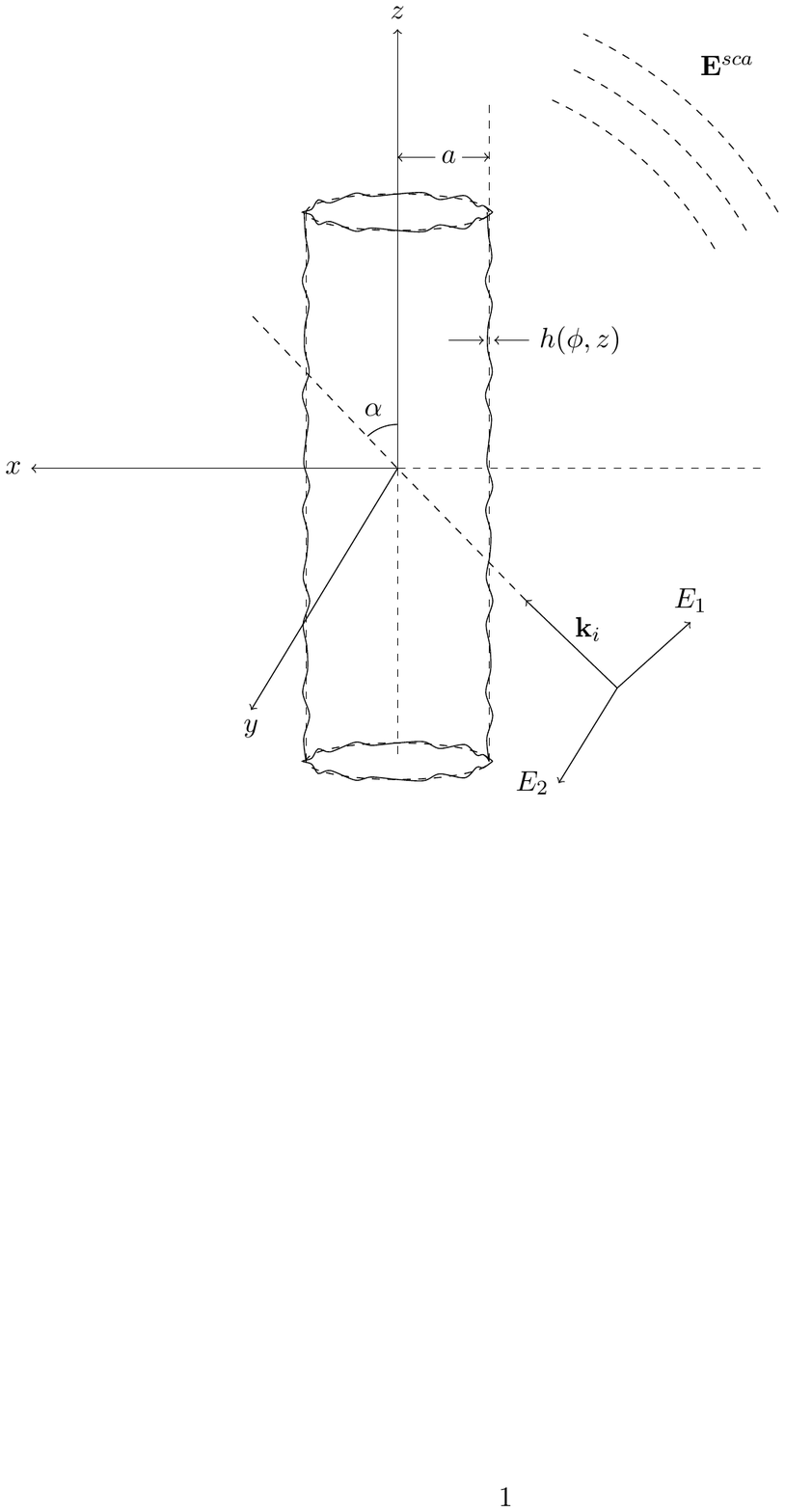}
\caption{Schematic of the rough cylindrical scatterer with a stochastic roughness $h(\phi,z)$ and illuminated by a plane wave propagating with wavevector $\textbf{k}_i$ in the $xz$ plane at an angle of $\alpha$ with the $x$ axis. $E_1$ is the component of the electric field in the $xz$ plane, while $E_2$ is the component along the $y$ axis.}
\label{scheme}
\end{figure}

Consider an azimuthally rough cylinder with it's surface being described by $r(\phi,z) = r_S(\phi)=a+h(\phi)$ (where $(r,\phi,z)$ are the cylindrical coordinates of the point in question), and as illustrated in the schematic in Fig.~ \ref{scheme}. The dielectric constant of the cylinder is assumed to be $\epsilon_d$ and it's permeability $\mu_d$. It is of interest to calculate the scattered fields when a plane wave is incident on this cylinder. We consider only TM polarisation ($\textbf{E} || \hat{z}$ and $\textbf{H} \perp \hat{z}$), which corresponds to an incident field with $\alpha = \pi/2$ and $E_2 = 0$ in Fig. ~\ref{scheme}. A harmonic time dependence of the form $\exp(j\omega t)$ is assumed in all the fields and suppressed throughout this tutorial. The incident field is given by:
\begin{equation}
\textbf{E}^{inc} = E_z^{inc} \hat{z}  = \exp(-jk_0(x\cos\phi_0 +y\sin \phi_0))\hat{z} = \exp(-jk_0 r \cos(\phi-\phi_0))\hat{z}
\end{equation}
which is a plane wave propagating at an angle $\phi_0$ with the $x$ axis and polarized in the $z$ direction. The starting point of our analysis is the electric field integral equation (EFIE) for homogeneous dielectric scatterers \cite[Ch.~1.9]{q}:
\begin{subequations}\label{EFIE}
\begin{align}
E_z^{inc} \bigg |_{S} = K(\phi) + j\eta_0 k_0 A_z^{(0)} \bigg |_S + \bigg[\frac{\partial F_y^{(0)}}{\partial x} - \frac{\partial F_x^{(0)}}{\partial y} \bigg]\bigg |_{S^{+}}\\
0 = -K(\phi) + j\eta_d k_d A_z^{(d)} \bigg |_S + \bigg[\frac{\partial F_y^{(d)}}{\partial x} - \frac{\partial F_x^{(d)}}{\partial y} \bigg]\bigg |_{S^{-}}
\end{align} 
\end{subequations}
where $\textbf{A}^{(p)} = A_z^{(p)} \hat{z}$ and $\textbf{F}^{(p)} = F_x^{(p)}\hat{x}+F_y^{(p)}\hat{y}$ for $p \in \{0,d\}$ are the magnetic and electric vector potential respectively defined by (the superscripts indicate the physical region in which the expression is valid: $0$ corresponds to vacuum and $d$ corresponds to dielectric):
\begin{subequations}\label{AF}
\begin{align}
&A_z^{(p)}(x,y) = \frac{1}{4j} \int \limits_0^{2\pi} J(\phi') H_0^{(2)}(k_p R) \sqrt{r_S^2(\phi')+r_S^{'2}(\phi')} \ d\phi' \\
&\textbf{F}^{(p)}(x,y) = \frac{1}{4j} \int \limits_0^{2\pi}\hat{t}(\phi') K(\phi') H_0^{(2)}(k_p R) \sqrt{r_S^2(\phi')+r_S^{'2}(\phi')} \ d\phi' \label{F}
\end{align}
\end{subequations}
where $J(\phi)$ and $K(\phi)$ are the equivalent electric and magnetic currents along the cylinder interface, and $R = [(x-r_S(\phi')\cos\phi')^2+(y-r_S(\phi')\sin\phi')^2]^{1/2}$ is the distance between a point $(r_S(\phi'),\phi')$ on the cylinder surface to the point of observation $(x,y)$. $\hat{t}(\phi)$ denotes a unit vector tangential to the cylinder surface at the point $(r_S(\phi),\phi)$ and is given by:
\begin{equation}
\hat{t}(\phi) = \frac{\hat{x}(-r_S(\phi)\sin \phi+r_S'(\phi)\cos \phi)+\hat{y}(r_S(\phi)\cos \phi+r_S'(\phi)\sin \phi)}{\sqrt{r_S^2(\phi)+r_S^{'2}(\phi)}}\label{unitvec}
\end{equation}
Using Eqs. (\ref{F}) and (\ref{unitvec}), it can be shown that:
\begin{align}
&\bigg[ \frac{\partial F_y^{(p)}}{\partial x} - \frac{\partial F_x^{(p)}}{\partial y}\bigg] \bigg |_{S} \nonumber \\&= \frac{k_p}{4j}\int \limits_0^{2\pi} K(\phi') \frac{H_0^{(2)'}(k_pR)}{R} (r_S(\phi')r_S(\phi)\cos(\phi-\phi')+r_S'(\phi')r_S(\phi)\sin(\phi'-\phi)-r_S^2(\phi')) \ d\phi'
\end{align}
To numerically implement a MoM routine to solve the EFIE (Eq.~(\ref{EFIE})), we assume a piecewise linear approximation of the cylinder surface (in $\phi$) between $N$ discrete points $(r_{n+1/2},\phi_{n+1/2})$, where $\phi_{n+1/2} = \frac{2\pi}{N}(n+\frac{1}{2})$ and $r_{n+1/2} = r_S(\phi_{n+1/2})$ $\forall \ n\in \{0,1,2,\dots N-1\}$. Thus, in this approximation, the radial distance of a point on the cylinder surface at an azimuthal angle $\phi$ satisfying $\phi_{n-1/2}\leq\phi<\phi_{n+1/2}$ is given by:
\begin{equation}\label{linear_inter}
r_S(\phi) \approx \frac{r_{n+1/2}+r_{n-1/2}}{2}+\frac{(r_{n+1/2}-r_{n-1/2})}{2\pi/N}\bigg (\phi-\frac{2\pi n}{N}\bigg) = r_n + s_n \bigg( \phi - \frac{2\pi n}{N}\bigg)
\end{equation}
where, for notational convenience, we have defined $r_n = (r_{n+1/2}+r_{n-1/2})/2$ and $s_n = \frac{N}{2\pi} (r_{n+1/2}-r_{n-1/2})$. In the analysis that follows, we denote the discretization angle $2\pi/N$ by $\delta$.

\begin{figure}[b]
\centering
\includegraphics[trim = 23cm 12.5cm 30cm 11cm, scale = 1.0]{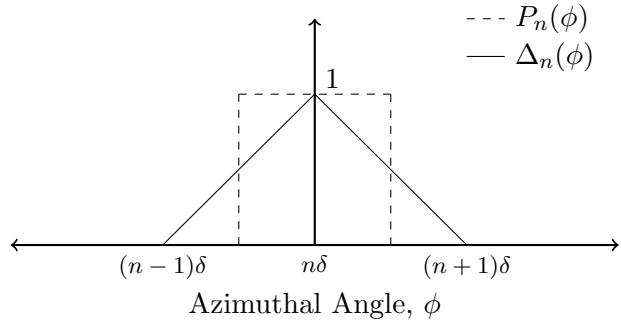}
\caption{Triangular basis functions, $\Delta_n(\phi)$ and pulse testing function, $P_n(\phi)$ used in the MoM implementation.}
\label{basis}
\end{figure}

The unknown functions $J(\phi)$ and $K(\phi)$ are then approximated as a linear combination of triangular basis functions (shown in Fig. \ref{basis}):
\begin{equation}\label{JK}
J(\phi) \approx \sum_{n=0}^{N-1} j_n \Delta_n(\phi) \ \text{and} \ K(\phi) \approx \sum_{n=0}^{N-1} k_n \Delta_n(\phi)
\end{equation}

where $ \forall n \in \{0,1,2\dots N-1\}$:
\begin{equation}
\Delta_n(\phi) = \begin{cases} 1-\frac{1}{\delta}\big|\phi-n\delta \big| \ &\text{for} \ (n-1)\delta<\phi<(n+1)\delta \\ 0 \ &\text{otherwise} \end{cases}
\end{equation}
Further, to impose Eq. (\ref{EFIE}), we employ the pulse testing functions $P_n(\phi)$ (see Fig. \ref{basis}) defined by ($\forall n\in\{0,1,2, \dots N-1\}$):
\begin{equation}
P_n(\phi) = \begin{cases}1 \ &\text{for} \ \big(n-\frac{1}{2}\big)\delta<\phi<\big(n+\frac{1}{2}\big)\delta \\ 0 \ &\text{otherwise}   \end{cases}
\end{equation}
Consider now the inner products of the individual terms in Eq. (\ref{EFIE}) with $P_m(\phi)$. For instance:
\begin{align}\label{amn}
\int\limits_0^{2\pi} A_z^{(p)} \bigg |_S P_m(\phi) \ d\phi =\frac{1}{4j} \int\limits_0^{2\pi} \int\limits_0^{2\pi} J(\phi')P_m(\phi)  H_0^{(2)}(k_pR){\sqrt{r_S^2(\phi')+r_S^{'2}(\phi')}} \ d\phi \ d\phi' \nonumber \\
= \sum_{n=0}^{N-1} \frac{j_n}{4j}  \int\limits_0^{2\pi} \int\limits_0^{2\pi} \Delta_n(\phi')P_m(\phi)  H_0^{(2)}(k_pR){\sqrt{r_S^2(\phi')+r_S^{'2}(\phi')}} \ d\phi \ d\phi' \triangleq \sum_{m=0}^{N-1}\alpha^{(p)}_{m,n}j_n
\end{align}
where
\begin{equation}\label{alpha}
\alpha_{m,n}^{(p)} =  \frac{1}{4j}  \int\limits_{(m-1/2)\delta}^{(m+1/2)\delta} \int\limits_{(n-1)\delta}^{(n+1)\delta} \Delta_n(\phi')H_0^{(2)}(k_pR){\sqrt{r_S^2(\phi')+r_S^{'2}(\phi')}} \ d\phi \ d\phi' 
\end{equation}
wherein we have used the fact that $P_m(\phi)$ is supported only on $((m-1/2)\delta,(m+1/2)\delta)$ and $\Delta_n(\phi)$ is supported on $((n-1)\delta,(n+1)\delta)$.
It is usually not possible to obtain a closed form expression for the integral in Eq. (\ref{alpha}). However, for sufficiently small $\delta$, the integral $\alpha_{m,n}^{(p)}$ can be approximated as:
\begin{equation}\label{break}
\alpha_{m,n}^{(p)} \approx \frac{1}{4j} \sqrt{r_{n-3/4}^2+s_{n-1}^2}\alpha_{m,n,1}^{(p)} + \frac{1}{4j} \sqrt{r_{n}^2+s_{n}^2}\alpha_{m,n,2}^{(p)}+\frac{1}{4j} \sqrt{r_{n+3/4}^2+s_{n+1}^2}\alpha_{m,n,3}^{(p)}
\end{equation}
with
\begin{subequations}\label{parts}
\begin{align}
\alpha_{m,n,1}^{(p)} =   \int\limits_{(m-1/2)\delta}^{(m+1/2)\delta} \int\limits_{(n-1)\delta}^{(n-1/2)\delta} \Delta_n(\phi') H_0^{(2)}(k_pR) \ d\phi \ d\phi' \\
\alpha_{m,n,2}^{(p)} =   \int\limits_{(m-1/2)\delta}^{(m+1/2)\delta} \int\limits_{(n-1/2)\delta}^{(n+1/2)\delta} \Delta_n(\phi') H_0^{(2)}(k_pR)\ d\phi \ d\phi' \\
\alpha_{m,n,3}^{(p)} =    \int\limits_{(m-1/2)\delta}^{(m+1/2)\delta} \int\limits_{(n+1/2)\delta}^{(n+1)\delta} \Delta_n(\phi') H_0^{(2)}(k_pR) \ d\phi \ d\phi'
\end{align}
\end{subequations}
where $r_{n-3/4} = (3r_{n-1/2}+r_{n-3/2})/4, r_{n+3/4} = (3r_{n+1/2}+r_{n+3/2})/4$. In writing Eq.~(\ref{break}), we have split the integral over $\phi'$ into three different integrals and approximated $\sqrt{r_S^2(\phi')+r_S^{'2}(\phi')}$ centroidally in each of them. Depending upon the values of $m$ and $n$, the integrals may or may not be singular. We deal with them in four cases:
\begin{enumerate}
\item For $m\neq n, n-1, n+1$, the integrands in Eq. (\ref{parts}) are not singular, and for sufficiently small $\delta$, they can be approximated centroidally:
\begin{subequations}
\begin{align}
\alpha_{m,n,1}^{(p)} =   \frac{\delta^2}{8} H_0^{(2)}(k_p R_{m,n-3/4})  \\
\alpha_{m,n,2}^{(p)} =  \frac{3\delta^2}{4} H_0^{(2)}(k_p R_{m,n})\\
\alpha_{m,n,3}^{(p)} =   \frac{\delta^2}{8} H_0^{(2)}(k_p R_{m,n+3/4})
\end{align}
\end{subequations}

where $R_{m,n} = \sqrt{r_m^2+r_n^2-2r_mr_n \cos((m-n)\delta)}$.

\item For $m = n$, the integrands in Eq. (\ref{parts}) are singular and have to be integrated analytically. To do so, we use the small argument approximation to the Hankel functions:
\begin{equation}\label{small_arg}
H_0^{(2)}(x) \approx 1-\frac{2j}{\pi} \ln\bigg(\frac{\gamma x}{2} \bigg) \ \text{for} \ x \to 0
\end{equation}
where $\gamma \approx 1.78107$ is the Euler's constant. Consider the integral in $\alpha_{m,m,2}^{(p)}$. Using Eq. (\ref{linear_inter}), $R \approx (\sqrt{r_m^2+s_m^2})|\phi-\phi'|$. Together with the small argument approximation, the integral becomes:
\begin{align}
\alpha_{m,m,2}^{(p)} &=    \int\limits_{(m-1/2)\delta}^{(m+1/2)\delta} \int\limits_{(m-1/2)\delta}^{(m+1/2)\delta} \Delta_n(\phi') \bigg[ 1-\frac{2j}{\pi} \ln\bigg(\frac{\gamma k_p \sqrt{r_m^2+s_m^2}|\phi-\phi'|}{2} \bigg)\bigg] \ d\phi \ d\phi' \nonumber \\
&= \frac{3\delta^2}{4} \bigg[ 1-\frac{2j}{\pi} \ln \bigg(\frac{\gamma k_p}{2}\sqrt{r_m^2+s_m^2} \bigg) \bigg]-\frac{2j\delta^2}{\pi} \int\limits_{-\delta/2}^{\delta/2} \int \limits_{-\delta/2}^{\delta/2} \ln |\phi-\phi'| \Delta_0(\phi')\ d\phi\ d\phi' \label{part}
\end{align}
The integral in Eq. (\ref{part}) can be evaluated analytically to give:
\begin{equation}
\alpha_{m,m,2}^{(p)} = \frac{3\delta^2}{4} \bigg[ 1-\frac{2j}{\pi} \ln \bigg(\frac{\gamma k_p}{2}\sqrt{r_m^2+s_m^2} \bigg) \bigg]-\frac{2j\delta^2}{\pi} \bigg[ \frac{3}{4}\ln(\delta) - 1.15057\bigg] \end{equation}
A similar procedure can be adopted to evaluate $\alpha_{m,m,1}^{(p)}$ and $\alpha_{m,m,3}^{(p)}$. We give the final result below:
\begin{subequations}
\begin{align}
\alpha_{m,m,1}^{(p)} =\frac{\delta^2}{8}\bigg[ 1-\frac{2j}{\pi} \ln \bigg(\frac{\gamma k_p}{2} \sqrt{r_{m-1/2}^2 +\frac{1}{4}(s_m+s_{m-1})^2} \bigg)\bigg]-\frac{2j\delta^2}{\pi}\bigg[\frac{1}{8}\ln(\delta) - 0.0699 \bigg]\\
\alpha_{m,m,3}^{(p)} = \frac{\delta^2}{8} \bigg[ 1-\frac{2j}{\pi} \ln \bigg(\frac{\gamma k_p}{2} \sqrt{r_{m+1/2}^2 +\frac{1}{4}(s_m+s_{m+1})^2} \bigg)\bigg]-\frac{2j\delta^2}{\pi} \bigg[\frac{1}{8}\ln(\delta) - 0.0699 \bigg]
\end{align}
\end{subequations}
\item For $m = n+1$, $\alpha_{m,m-1,1}^{(p)}$ can be evaluated using the centroidal approximation, since the integrand is not singular:
\begin{equation}
\alpha_{m,m-1,1}^{(p)} = \frac{\delta^2}{8} H_0^{(2)}(k_p R_{m,m-7/4})
\end{equation}
A procedure similar to that used in evaluating $\alpha_{m,m,2}^{(p)}$ can be used to evaluate $\alpha_{m,m-1,2}$ and $\alpha_{m,m-1,3}$:
\begin{subequations}
\begin{align}
&\alpha_{m,m-1,2}^{(p)} = \frac{3\delta^2}{4}\bigg[1-\frac{2j}{\pi}\ln\bigg( \frac{\gamma k_p}{2} \sqrt{r_{m-1/2}^2+\frac{1}{4}(s_m+s_{m-1})^2} \bigg) \bigg]-\frac{2j\delta^2}{\pi} (0.749 \ln (\delta)+0.33)\\
&\alpha_{m,m-1,3}^{(p)} = \frac{\delta^2}{8}\bigg[1-\frac{2j}{\pi} \ln\bigg( \frac{\gamma k_p}{2} \sqrt{r_{m}^2+s_m^2} \bigg) \bigg]-\frac{2j\delta^2}{\pi} (0.3124 \ln (\delta)-0.29024)
\end{align}
\end{subequations}
\item The case of $m = n-1$ is identical to the previous case:
\begin{subequations}
\begin{align}
&\alpha_{m,m+1,1}^{(p)} = \frac{\delta^2}{8}\bigg[1-\frac{2j}{\pi} \ln\bigg( \frac{\gamma k_p}{2} \sqrt{r_{m}^2+s_m^2} \bigg) \bigg]-\frac{2j\delta^2}{\pi} (0.3124 \ln (\delta)-0.29024)\\
&\alpha_{m,m+1,2}^{(p)} = \frac{3\delta^2}{4}\bigg[1-\frac{2j}{\pi}\ln\bigg( \frac{\gamma k_p}{2} \sqrt{r_{m+1/2}^2+\frac{1}{4}(s_m+s_{m+1})^2} \bigg) \bigg]-\frac{2j\delta^2}{\pi} (0.749 \ln (\delta)+0.33)\\
&\alpha_{m,m+1,3}^{(p)} =\frac{\delta^2}{8} H_0^{(2)}(k_p R_{m,m+7/4})
\end{align}
\end{subequations}
\end{enumerate}
Consider now the inner product of $\bigg[ \frac{\partial F_y^{(p)}}{\partial x} - \frac{\partial F_x^{(p)}}{\partial y}\bigg] \bigg |_{S}$ with $P_m(\phi)$:
\begin{align}\label{ban}
\int\limits_{0}^{2\pi} \bigg[ \frac{\partial F_y^{(p)}}{\partial x} - \frac{\partial F_x^{(p)}}{\partial y}\bigg] \bigg |_{S} P_m(\phi) \ d\phi \triangleq \sum_{n=0}^{N-1} \beta_{m,n}^{(p)}k_n
\end{align}
where
\begin{equation}
\beta_{m,n}^{(p)} =  \frac{k_p}{4j}\int \limits_{(m-1/2)\delta}^{(m+1/2)\delta} \int \limits_{(n-1)\delta}^{(n+1)\delta} \Delta_{n}(\phi') \frac{H_0^{(2)'}(k_pR)}{R} (r_S(\phi')r_S(\phi)\cos(\phi-\phi')+r_S'(\phi')r_S(\phi)\sin(\phi'-\phi)-r_S^2(\phi')) \ d\phi'\ d \phi
\end{equation}
which we again express as a sum of three integrals: $\beta_{m,n}^{(p)} = \frac{k_p}{4j} (\beta_{m,n,1}^{(p)}+\beta_{m,n,2}^{(p)}+\beta_{m,n,3}^{(p)})$ where:
\begin{subequations}\label{partb}
\begin{align}
\beta_{m,n,1}^{(p)} = \int \limits_{(m-1/2)\delta}^{(m+1/2)\delta} \int \limits_{(n+1/2)\delta}^{(n+1)\delta} \Delta_{n}(\phi') \frac{H_0^{(2)'}(k_pR)}{R} F(\phi,\phi') \ d\phi' \ d\phi\\
\beta_{m,n,2}^{(p)} = \int \limits_{(m-1/2)\delta}^{(m+1/2)\delta} \int \limits_{(n-1/2)\delta}^{(n+1/2)\delta} \Delta_{n}(\phi') \frac{H_0^{(2)'}(k_pR)}{R} F(\phi,\phi') \ d\phi' \ d\phi\\
\beta_{m,n,3}^{(p)} = \int \limits_{(m-1/2)\delta}^{(m+1/2)\delta} \int \limits_{(n-1)\delta}^{(n-1/2)\delta} \Delta_{n}(\phi') \frac{H_0^{(2)'}(k_pR)}{R} F(\phi,\phi') \ d\phi' \ d\phi
\end{align}
\end{subequations}
where $F(\phi,\phi') = r_S(\phi')r_S(\phi)\cos(\phi-\phi')+r_S'(\phi')r_S(\phi)\sin(\phi-\phi')-r_S^2(\phi')$. We again consider four different cases:
\begin{enumerate}
\item For $m \neq n, n-1,n+1$, the integrals in Eq.~(\ref{partb}) are not singular and can be approximated centroidally:
\begin{align}
&\beta_{m,n,1}^{(p)} =\frac{\delta^2}{8} \frac{H_0^{(2)'}(k_p R_{m,n+3/4})}{R_{m,n+3/4}} F_{m,n+3/4} \\
&\beta_{m,n,2}^{(p)} =\frac{3\delta^2}{4} \frac{H_0^{(2)'}(k_p R_{m,n})}{R_{m,n}} F_{m,n} \\
&\beta_{m,n,3}^{(p)} =\frac{\delta^2}{8} \frac{H_0^{(2)'}(k_p R_{m,n-3/4})}{R_{m,n-3/4}} F_{m,n-3/4}
\end{align}
where $F_{m,n} = r_{m}r_n \cos((n-m)\delta)+r_m s_n \sin((m-n)\delta)-r_n^2$ with $s_n$ denoting the slope of the segment containing $r_n$. For instance $s_{n-3/4} = (r_{n-1/2}-r_{n-3/2})/\delta$.
\item For $m = n$, the integral is singular. Additionally, note from Eq.~(\ref{EFIE}) that what needs to be evaluated is not $[ \frac{\partial F_y^{(p)}}{\partial x} - \frac{\partial F_x^{(p)}}{\partial y}] |_{S}$, rather $[ \frac{\partial F_y^{(0)}}{\partial x} - \frac{\partial F_x^{(0)}}{\partial y}] |_{S^+}$ or $[ \frac{\partial F_y^{(d)}}{\partial x} - \frac{\partial F_x^{(d)}}{\partial y}] |_{S^-}$. If the integrand is not singular, then the limit is the same as the value of the inner product at the surface $S$, however care must be taken in computing the limit for the case when the integrand is singular. Consider, for instance, the evaluation of $\beta_{m,m,2}^{(0)}$:
\begin{equation}
\beta_{m,m,2}^{(0)} = \lim_{\epsilon\to 0} \int \limits_{(m-1/2)\delta}^{(m+1/2)\delta} \int \limits_{(m-1/2)\delta}^{(m+1/2)\delta} \Delta_{n}(\phi') \frac{H_0^{(2)'}(k_0R)}{R} F(\phi,\phi') \ d\phi' \ d\phi
\end{equation}
with
\begin{align}
R = \sqrt{r_s(\phi')^2+(r_S(\phi)+\epsilon)^2-r_S(\phi')(r_S(\phi)+\epsilon) \cos(\phi-\phi')} 
 \approx \sqrt{r_m^2(\phi-\phi')^2+\epsilon^2}
\end{align}
wherein we have used the fact that since $|\phi-\phi'|<\delta$, $|\phi-\phi'|$ is small and thus $\cos(\phi-\phi') \approx 1-\frac{(\phi-\phi')^2}{2}$ and $\sin(\phi-\phi') \approx (\phi-\phi')$.   Also note that $\epsilon>0$, since $\beta_{m,m,2}^{(0)}$ arises in the evaluation of $[ \frac{\partial F_y^{(0)}}{\partial x} - \frac{\partial F_x^{(0)}}{\partial y}] |_{S^+}$. We may also use the small argument approximation and set $H_0^{(2)'}(k_0R) \approx -\frac{2j}{\pi k_0 R}$. Moreover, for small $|\phi-\phi'|$, the term $F(\phi,\phi')$ can be approximated by $[\epsilon r_m - \frac{r_m^2}{2}(\phi-\phi')^2]$. Therefore:
\begin{align}
\beta_{m,m,2}^{(0)} &= \frac{j}{\pi k_0}  \lim_{\epsilon\to 0}\bigg[ \int \limits_{(m-1/2)\delta}^{(m+1/2)\delta} \int \limits_{(m-1/2)\delta}^{(m+1/2)\delta} \bigg(1+\bigg|\frac{\phi'-m \delta}{\delta}\bigg| \bigg)\frac{r_m^2(\phi-\phi')^2-2\epsilon r_m}{r_m^2 (\phi-\phi')^2+\epsilon^2} \ d\phi' \ d\phi \bigg] \nonumber \\ &= \frac{j}{\pi k_0}  \lim_{\epsilon\to 0}\bigg[ \int \limits_{-\delta/2}^{\delta/2} \int \limits_{-\delta/2}^{\delta/2} \bigg(1+\bigg|\frac{\phi'}{\delta}\bigg| \bigg)\frac{r_m^2(\phi-\phi')^2-2\epsilon r_m}{r_m^2 (\phi-\phi')^2+\epsilon^2} \ d\phi' \ d\phi \bigg]\label{integral}
 \end{align}
where we have performed the substitution $\phi - m\delta \to \phi$ and $\phi'-m\delta \to \phi'$ in the last step. Consider the integral:
\begin{align}
 \lim_{\epsilon\to 0}  \int \limits_{-\delta/2}^{\delta/2}\frac{r_m^2(\phi-\phi')^2-2\epsilon r_m}{r_m^2 (\phi-\phi')^2+\epsilon^2}  d \phi &= \delta -2\lim_{\epsilon \to 0} \int\limits_{-\delta/2}^{\delta/2} \frac{\epsilon/r_m}{(\phi-\phi')^2+\epsilon^2/r_m^2} d \phi \nonumber \\
& = \delta - 2 \lim_{\epsilon \to 0}\bigg[\tan^{-1}\bigg(\frac{r_m}{\epsilon}\bigg(\frac{\delta}{2}-\phi' \bigg)\bigg)+\tan^{-1}\bigg(\frac{r_m}{\epsilon}\bigg(\frac{\delta}{2}+\phi' \bigg)\bigg) \bigg]
\end{align}
Since, in Eq.~(\ref{integral}), $|\phi'|<\delta/2$,
\begin{equation}
 \lim_{\epsilon\to 0}  \int \limits_{-\delta/2}^{\delta/2}\frac{r_m^2(\phi-\phi')^2-2\epsilon r_m}{r_m^2 (\phi-\phi')^2+\epsilon^2}  d \phi =  \delta - 2\pi
\end{equation}
Note the extra $2\pi$ term in the above integral. Had we taken the limit before evaluating the integral, we would have obtained $\delta$ as the result of the above integration, which would yield erroneous results. Substituting this integral into Eq. (\ref{integral}), we obtain:
\begin{equation}
\beta_{m,m,2}^{(0)} = \frac{3j}{4\pi k_0} \delta(\delta-2\pi)
\end{equation}
Working similarly with $\beta_{m,m,2}^{(d)}$ (where it would be necessary to use $\epsilon<0$ since $\beta_{m,m,2}^{(d)}$ arises in the evaluation of $[ \frac{\partial F_y^{(d)}}{\partial x} - \frac{\partial F_x^{(d)}}{\partial y}] |_{S^-}$), we would obtain:
\begin{equation}
\beta_{m,m,2}^{(d)} = \frac{3j}{4\pi k_d} \delta(\delta+2\pi)
\end{equation}
$\beta_{m,m,1}^{(p)}$ and $\beta_{m,m,3}^{(p)}$ can be evaluated following the same procedure, resulting in:
\begin{equation}
\beta_{m,m,1}^{(p)} = \beta_{m,m,3}^{(p)} = \frac{j \delta^2}{8\pi k_p}
\end{equation}
Note that in the above mentioned calculations, we have ignored the slope $s_n$ of the surface. For a sufficiently fine sampling, this approximation is a valid one, since $r_{n-1/2} \approx r_{n+1/2}$, thereby implying $s_n \approx 0$.
\item For $n = m+1$, the integrand in $\beta_{m,m+1,1}^{(p)}$ is not singular and can be approximated centroidally:
\begin{equation}
\beta_{m,m+1,1}^{(p)} = \frac{\delta^2}{8}\frac{H_0^{(2)'}(k_p R_{m,m+7/4})}{R_{m,m+7/4}} F_{m+7/4,m}
\end{equation}
$\beta_{m,m+1,2}^{(p)}$ is singular and must be evaluated using the procedure described above. It can be shown that in this case, the limiting value of the integral is same as the value of the integral at the cylinder surface:
\begin{equation}
\beta_{m,m+1,2}^{(p)} = \frac{3j\delta^2}{4\pi k_p}
\end{equation}
Care must be taken in evaluating the limit in $\beta_{m,m+1,3}^{(p)}$:
\begin{align}
\beta_{m,m+1,3}^{(p)} = \begin{cases}\frac{j \delta}{8\pi k_0} (\delta-2\pi) &\text{for} \ p = 0 \\ \frac{j \delta}{8\pi k_d} (\delta+2\pi) &\text{for} \ p = d  \end{cases}
\end{align}
\item The case of $m = n+1$ is identical to the case in which $m = n-1$. We give the final results below:
\begin{subequations}
\begin{align}
&\beta_{m,m-1,1}^{(p)} =  \begin{cases}\frac{j \delta}{8\pi k_0} (\delta-2\pi) &\text{for} \ p = 0 \\ \frac{j \delta}{8\pi k_d} (\delta+2\pi) &\text{for} \ p = d  \end{cases} \\
&\beta_{m,m-1,2}^{(p)} = \frac{3j\delta^2}{4\pi k_p} \\
&\beta_{m,m-1,3}^{(p)} = \frac{\delta^2}{8}\frac{H_0^{(2)'}(k_p R_{m,m-7/4})}{R_{m,m-7/4}} F_{m-7/4,m}
\end{align}
\end{subequations}
\end{enumerate}
Consider the inner product of $E_z^{(inc)} |_S$ with $P_m(\phi)$, which can be approximated centroidally:
\begin{equation}\label{en}
\int \limits_0^{2\pi} E_z^{(inc)} \bigg |_S P_m(\phi) \ d\phi \triangleq e_m \approx \delta \exp(-jk_0 r_m \cos(m\delta-\phi_0))
\end{equation}
Finally consider the inner product of $K(\phi)$ with $P_m(\phi)$:
\begin{align}
\int\limits_0^{2\pi} K(\phi) P_m(\phi) \ d\phi = \sum_{n=0}^{N-1} k_n \int\limits_{(m-1/2)\delta}^{(m+1/2)\delta} \Delta_n(\phi) P_m(\phi) = \frac{3\delta}{4}k_m +\frac{\delta}{8}(k_{m-1}+k_{m+1})
\end{align}
Imposing Eq.~(\ref{EFIE}) by taking the inner product of LHS and RHS with $P_m(\phi)$ for $m \in \{0,1,2,\dots N-1\}$, we obtain the following system of linear equations:
\begin{subequations}\label{final}
\begin{align}
e_m = \frac{3\delta}{4}k_m +\frac{\delta}{8}(k_{m-1}+k_{m+1})+j\eta_0 k_0 \sum_{n=0}^{N-1} \alpha_{nm}^{(0)} j_m +\sum_{n=0}^{N-1}\beta_{nm}^{(0)} k_m \\
0 = -\frac{3\delta}{4}k_m +\frac{\delta}{8}(k_{m-1}+k_{m+1})+j\eta_d k_d \sum_{n=0}^{N-1} \alpha_{nm}^{(d)} j_m +\sum_{n=0}^{N-1}\beta_{nm}^{(d)} k_m 
\end{align}
\end{subequations}
which is a system of $2N$ linear equations in $2N$ unknowns ($j_n , k_n \ \forall \ n\in\{0,1,2,3 \dots N-1\}$) that can be solved to obtain $j_n$ and $k_n$ and thus $J(\phi)$ and $K(\phi)$ through Eq. (\ref{JK}). Once $J(\phi)$ and $K(\phi)$ are known, the scattered fields can be evaluated through:
\begin{equation}\label{scaf}
E_z^{sca} = -j\eta_0 k_0A_z^{(0)} - \bigg[ \frac{\partial F_y^{(0)}}{\partial x} - \frac{\partial F_x^{(0)}}{\partial y}\bigg]
\end{equation}
wherein $A_z$ and $\textbf{F}$ are defined through Eq. (\ref{AF}). A quantity of interest related to the scatterer is the two-dimensional scattering cross-section $\sigma(\phi)$ defined by:
\begin{equation}
\sigma(\phi) = \lim_{r\to \infty} 2\pi r \bigg| \frac{(E_z^{sca})^2}{E_0^2}\bigg|
\end{equation}
$E_0$ being the amplitude of the incident plane wave. The limiting form of $E_z^{sca}$ as $r\to \infty$ can be evaluated using Eqs.~(\ref{scaf}) and (\ref{AF}). Note that as $r\to\infty$,
\begin{equation}
R = \sqrt{r^2+r_S(\phi')^2-2rr_S(\phi') \cos(\phi-\phi')} \approx r - r_S(\phi')\cos(\phi-\phi')
\end{equation}
and thus
\begin{subequations}
\begin{align}
&H_0^{(2)}(k_0 R) \approx \sqrt{\frac{2}{\pi k_0 r}}\exp(-j(k_0 r-\pi/4)) \exp(jk_0 r_S(\phi')\cos(\phi-\phi')) \\ 
&H_0^{(2)'}(k_0 R) \approx -j\sqrt{\frac{2}{\pi k_0 r}}\exp(-j(k_0 r-\pi/4)) \exp(jk_0 r_S(\phi')\cos(\phi-\phi')) 
\end{align}
\end{subequations}
Similarly,
\begin{equation}
F(\phi,\phi') \approx r (r_S(\phi')\cos(\phi-\phi')+r_S'(\phi')\sin(\phi'-\phi))
\end{equation}
Using Eqs. (\ref{AF}) and (\ref{scaf}), we get:
\begin{align}
\lim _{r\to \infty} E_z^{sca} &= -\sqrt{\frac{2}{\pi k_0 r}} \frac{k_0}{4}\exp(-jk_0(r-\pi/4)) \bigg[\int\limits_0^{2\pi} \sqrt{r_S^2(\phi')+r_S^{2'}(\phi')} \eta_0J(\phi')\nonumber \\ &-K(\phi') (r_S(\phi')\cos(\phi-\phi')+r_S'(\phi')\sin(\phi'-\phi))\exp(jk_0 r_S(\phi')\cos(\phi-\phi') \ d\phi' \bigg] 
\end{align}
and thus
\begin{align}
\sigma(\phi) &= \frac{k_0}{4 E_0^2}\bigg| \bigg[\int\limits_0^{2\pi} \sqrt{r_S^2(\phi')+r_S^{2'}(\phi')} \eta_0J(\phi')-K(\phi') (r_S(\phi')\cos(\phi-\phi')+r_S'(\phi')\sin(\phi'-\phi))\nonumber \\ &\exp(jk_0 r_S(\phi')\cos(\phi-\phi') \ d\phi' \bigg]  \bigg|^2
\end{align}

 \begin{figure}[htpb]
 \centering
 \begin{subfigure}{\textwidth}
  \centering
\includegraphics[scale = 0.3]{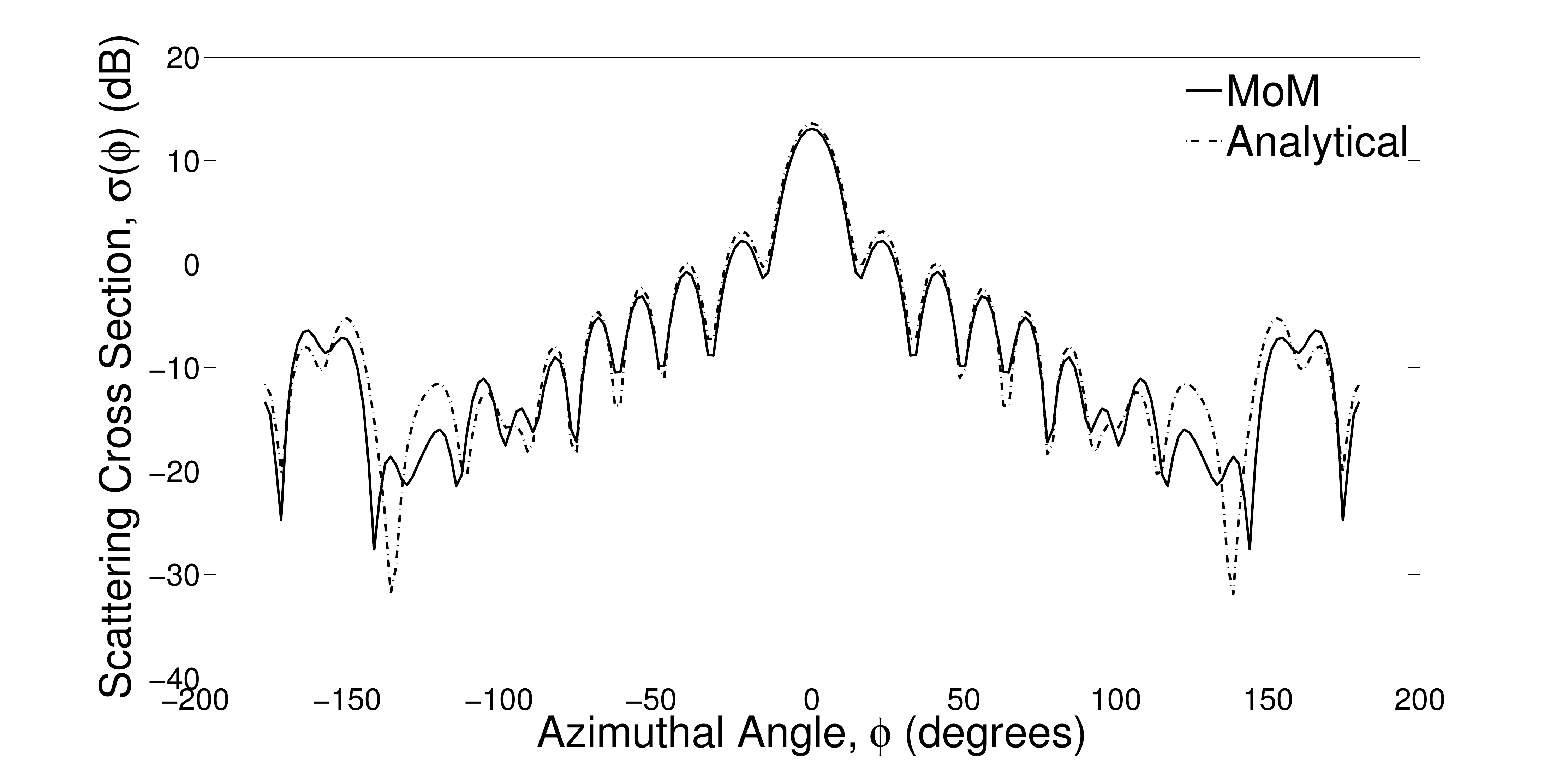}
\caption{$N = 100, \delta = 0.063$ rad}
\end{subfigure}
\begin{subfigure}{\textwidth}
 \centering
\includegraphics[scale = 0.3]{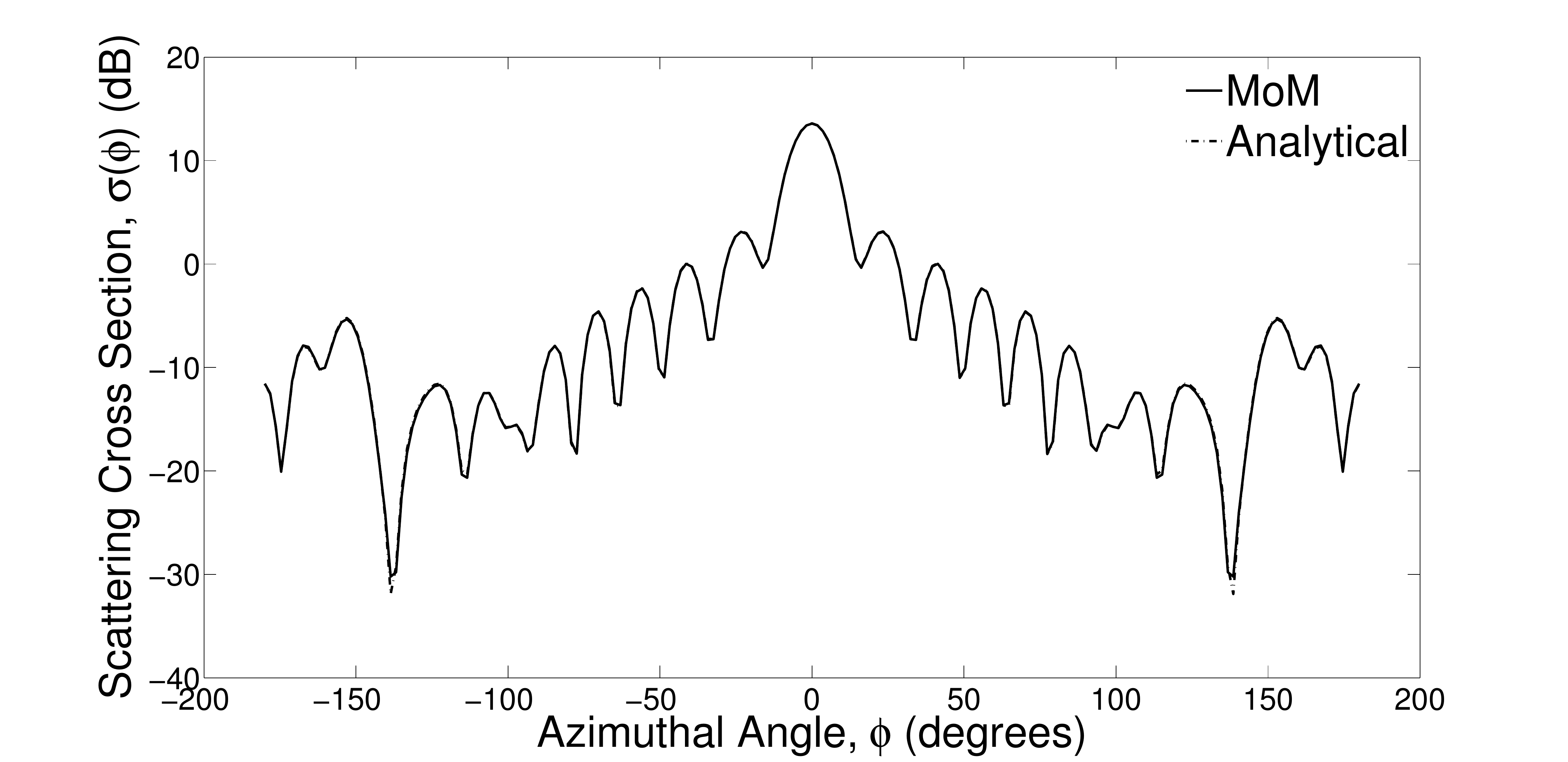}
\caption{$N = 300, \delta = 0.021$ rad}
\end{subfigure}
\caption{Plot of $\sigma(\phi)$ computed using MoM and compared with the analytical solution for different discretizations of the cylinder surface. In all the calculations, $a =2\lambda_0$, permittivity $\epsilon_d = 2\epsilon_0$, and permeability $\mu_d = \mu_0$.}
\label{Azvarnh0}
\end{figure}

Using Eq. (\ref{JK}), $\sigma(\phi)$ can be expressed in terms of $j_n$ and $k_n$:
\begin{equation}
\sigma(\phi) = \frac{k_0}{4 E_0^2} \bigg|\sum_{n=0}^{N-1} (a_n(\phi) j_n - b_n(\phi) k_n) \bigg|^2
\end{equation}
where
\begin{subequations}\label{sca_ab}
\begin{align}
&a_n(\phi) =\eta_0 \int \limits_0^{2\pi} \Delta_n(\phi')  \sqrt{r_S^2(\phi')+r_S^{2'}(\phi')} \exp(jk_0 r_S(\phi') \cos(\phi-\phi')) \ d\phi' \\
&b_n(\phi) =\int \limits_0^{2\pi} \Delta_n(\phi')  (r_S(\phi')\cos(\phi-\phi')+r_S'(\phi')\sin(\phi'-\phi)) \exp(jk_0 r_S(\phi') \cos(\phi-\phi')) \ d\phi' 
\end{align}
\end{subequations}
Since the integrals in Eq. (\ref{sca_ab}) are not singular, they can be approximated centroidally:
\begin{subequations}
\begin{align}
&a_n(\phi) \approx \eta_0 \bigg[\frac{\delta}{8}\sqrt{r_{n-3/4}^2+s_{n-1}^2}\exp(jk_0 r_{n-3/4}\cos((n-3/4\delta)-\phi) +\frac{3\delta}{4} \sqrt{r_n^2+s_n^2} \times \nonumber \\ &\exp(jk_0 r_n \cos(n\delta-\phi)+\frac{\delta}{8}\sqrt{r_{n+3/4}^2+s_{n+1}^2}\exp(jk_0 r_{n+3/4}\cos((n+3/4\delta)-\phi) \bigg]
\end{align}
\begin{align}
&b_n(\phi) \approx \frac{\delta^2}{8} (r_{n-3/4}\cos((n-3/4)\delta-\phi) + s_{n-1} \sin ((n-3/4)\delta-\phi))\times  \nonumber \\
&\exp(jk_0 r_{n-3/4}\cos((n-3/4)\delta-\phi))+ \frac{3\delta^2}{4} (r_{n}\cos(n\delta-\phi) + s_{n} \sin (n\delta-\phi))\exp(jk_0 r_{n}\cos(n\delta-\phi)) \nonumber\\ & +\frac{\delta^2}{8} (r_{n+3/4}\cos((n+3/4)\delta-\phi) + s_{n+1} \sin ((n+3/4)\delta-\phi))\exp(jk_0 r_{n+3/4}\cos((n+3/4)\delta-\phi))
\end{align}
\end{subequations}

Fig.~\ref{Azvarnh0} shows a comparison between the MoM results and analytical calculation based on the Mie theory \cite[Ch.~11.5]{balanis1989advanced} for a smooth cylinder $r_S(\phi) = a =$ const. It  can be seen that for $\delta \sim \frac{1}{25}\frac{\lambda_0}{a}$, the MoM described above can be used to accurately compute the scattered fields.

\section{Scattering from a perfectly conducting cylinder with axial roughness}\label{SecAxR}
Consider a perfectly conducting cylinder with it's axis oriented along the $z$ axis and the surface being described by $r(\phi,z) = r_S(z) =  a+h(z)$, where $h(z)$ is the roughness function and $a$ is the mean radius of the cylinder. We consider an azimuthally symmetric incident field given by:
\begin{equation}
\textbf{E}^{inc}(r,\phi,z) = E^{inc}(r,z)\hat{\phi}
\end{equation}
The above field corresponds to an incoming cylindrical wave incident on the cylinder. To compute the scattered fields, we use the EFIE at the surface of the PEC \cite[Ch.~1.7]{q}:
\begin{equation}\label{EFIE2}
\hat{n} \times \textbf{E}^{inc}\bigg |_S = \frac{j\eta_0}{k_0}\hat{n}\times(\nabla \nabla \cdot \textbf{A}+k_0^2 \textbf{A}) \bigg |_{S}
\end{equation}
where
\begin{equation}\label{def_A}
\textbf{A}(\textbf{r}) = \int \textbf{J}(\textbf{r}')\frac{\exp(-jk_0 |\textbf{r}-\textbf{r}'|)}{4\pi |\textbf{r}-\textbf{r}'|} d^3 \textbf{r}'
\end{equation} 
Symmetry dictates $\textbf{J}(\textbf{r}) = J(z) \hat{\phi}$. Let $\textbf{r} \equiv (r,\phi,z)$ and $\textbf{r}' \equiv (r',\phi',z')$. Noting that $\textbf{J}$ is actually a surface current, Eq.~\ref{def_A} can be written as ($R = [r^2+r_S(z')^2-2rr_S(z') \cos(\phi-\phi')+(z-z')^2]^{1/2}$):
\begin{align}
\textbf{A}(\textbf{r}) = \int\limits_{z'=-\infty}^{\infty} \int \limits_{\phi' = 0}^{2\pi} J(z')\hat{\phi}'\frac{\exp(-jk_0 R)}{4\pi R}\ r_S(z') \sqrt{1+r_S'^2(z')}d\phi' dz'
\end{align}
Observing that $\hat{\phi}' = -\sin \phi' \ \hat{x} + \cos\phi' \ \hat{y}$ and performing the substitution $\theta= \phi'-\phi$, we obtain (since the integrand is a periodic function of $\theta$ or $\phi'$, the exact limits are irrelevant and integration over any period will give the same result, thus the limits have been chosen from $-\pi$ to $\pi$ for convenience.):
\begin{align}\label{exp_A}
\textbf{A}(\textbf{r}) = \int\limits_{z'=-\infty}^{\infty} dz' J(z') r_S(z') \sqrt{1+r_S'^2(z')} \int \limits_{\phi' = -\pi}^{\pi} (-\sin( \phi+\theta) \hat{x} + \cos(\phi+\theta) \hat{y})\frac{\exp(-jk_0 R)}{4\pi R}\ d \theta
\end{align}
Also note that $R$ is an even function of $\theta$. Consider the integral over $\theta$:
\begin{align}\label{simp_A}
\int \limits_{\theta = -\pi}^{\pi} (-\sin( \phi+\theta) \hat{x} + \cos(\phi+\theta) \hat{y})\frac{\exp(-jk_0 R)}{4\pi R}\ d \theta = \hat{\phi} \int \limits_{\theta = -\pi}^{\pi} \cos\theta \frac{\exp(-jk_0 R)}{4\pi R}\ d \theta
\end{align}
with $\hat{\phi} = -\sin \phi\ \hat{x} + \cos\phi\ \hat{y}$. To arrive at this result, we expanded $\sin(\phi+\theta) = \sin \phi \cos\theta+\cos \phi \sin \theta$ and $\cos(\phi+\theta) = \cos \phi \cos\theta- \sin \phi \sin \theta$ and used the fact that:
\begin{equation}
 \int \limits_{\theta = -\pi}^{\pi} \sin\theta \frac{\exp(-jk_0 R)}{4\pi R}\ d \theta = 0
\end{equation}
since the integrand is an odd function of $\theta$. Plugging Eq. (\ref{simp_A}) into Eq. (\ref{exp_A}), and noting that $\hat{\phi}$ is a constant with respect to the integral, we obtain:

\begin{equation}\label{final_A}
\textbf{A}(\textbf{r}) \triangleq A(r,z) \hat{\phi} = \hat{\phi} \int\limits_{z'=-\infty}^{\infty} dz' J(z') r_S(z') \sqrt{1+r_S'^2(z')} \int \limits_{\theta = -\pi}^{\pi} \cos\theta \frac{\exp(-jk_0 R)}{4\pi R}\ d \theta
\end{equation}

Note that Eq. (\ref{final_A}) immediately implies that $\nabla \cdot \textbf{A} = 0$. Therefore, the EFIE (Eq. (\ref{EFIE2})) simplifies to:
\begin{equation}\label{red_EFIE}
E^{inc}(r,z)\bigg |_S = j\eta_0 k_0 A(r,z) \bigg |_S
\end{equation}
To numerically solve Eq. (\ref{red_EFIE}), we consider only  a finite length of the cylinder, say from $z = -L/2$ to $z = L/2$. As in Section \ref{SecAzR}, a piecewise linear approximation to the cylinder surface between $N+1$ discrete points is assumed. Let the points be labeled as $(z_n,r_n)$, where $n \in \{0,1,2, \dots N\}$ and $z_n = n\delta, r_n = r_S(z_n)$, with $\delta$ being the discretisation length $\delta = L/N$. The radial coordinate of a point on the cylinder surface at $z$, where $z_n < z<z_{n+1}$ is approximated by:
\begin{equation}
r_S(z) \approx \frac{z_{n+1}+z_{n}}{2} + \frac{z_{n+1} - z_{n}}{\delta} \bigg( z - \bigg(n+\frac{1}{2} \bigg)\delta\bigg) = z_{n+1/2}+s_n \bigg( z - \bigg(n+\frac{1}{2} \bigg)\delta\bigg)
\end{equation}
where we have defined $z_{n+1/2} = \frac{z_{n+1}+z_{n}}{2} $ and $s_n = \frac{z_{n+1} - z_{n}}{\delta} $. We expand the unknown $J(z)$ on the pulse basis functions:
\begin{equation}\label{J_appx}
J(z) \approx \sum_{n = 0}^{N-1} j_n P_{n+1/2}(z)
\end{equation}
and test Eq. (\ref{red_EFIE}) through the delta testing functions $T_{n+1/2}(z) = \delta(z-z_{n+1/2})$. Consider, now the inner product of each term in Eq. (\ref{red_EFIE}) with $T_{n+1/2}(z)$. For instance:
\begin{equation}
 \int \limits_0^{L} E^{inc}(r,z) \bigg |_S T_{m+1/2}(z) \ dz = E^{inc}(r_{m+1/2},z_{m+1/2}) \triangleq e_m
\end{equation}
where $r_{m+1/2} = r_S(z_{m+1/2}) \approx \frac{r_m+r_{m+1}}{2}$. Similarly:
\begin{align}
 \int \limits_0^{L} A(r,z) \bigg |_S T_{m+1/2}(z) \ dz = \int\limits_{z'=-\infty}^{\infty} dz' J(z') r_S(z') \sqrt{1+r_S'^2(z')} \int \limits_{\theta = -\pi}^{\pi} \cos\theta \frac{\exp(-jk_0 R_m)}{4\pi R_m}\ d \theta
 \triangleq \sum_{n=0}^{N-1} \alpha_{m,n} j_n
\end{align}
where
\begin{align}
\alpha_{m,n} &= \int\limits_{z'=-\infty}^{\infty} dz' P_{n+1/2}(z') r_S(z') \sqrt{1+r_S'^2(z')} \int \limits_{\theta = -\pi}^{\pi} \cos\theta \frac{\exp(-jk_0 R_m)}{4\pi R_m}\ d \theta \nonumber \\
 &= \int\limits_{z'=n\delta}^{(n+1)\delta} dz' r_S(z') \sqrt{1+r_S'^2(z')} \int \limits_{\theta = -\pi}^{\pi} \cos\theta \frac{\exp(-jk_0 R_m)}{4\pi R_m}\ d \theta 
\end{align}
with $R_m = R = [r_S(z')^2+r_{m+1/2}^2-2r_S(z')r_{m+1/2} \cos(\phi-\phi')+(z'-z_{m+1/2})^2]^{1/2}$. We treat the evaluation of $\alpha_{m,n}$ seperately for two cases:
\begin{enumerate}
\item For $m\neq n$, the integral in $z'$ can be approximated centroidally resulting in:
\begin{equation}
\alpha_{m,n} = \delta r_{n+1/2}\sqrt{1+s_n^2} \int \limits_{\theta = -\pi}^{\pi} \cos\theta \frac{\exp(-jk_0 R_{m,n})}{4\pi R_{m,n}}\ d \theta 
\end{equation}
where $R_{m,n} = [r_{m+1/2}^2+r_{n+1/2}^2+(m-n)^2\delta^2-2r_{n+1/2}r_{m+1/2}\cos \theta]^{1/2}$. The integral over $\theta$ in $\alpha_{m,n}$ would have to be evaluated numerically using quadrature methods.
\item For $m = n$, the integral is singular and the singularity would have to be integrated out. Following the method described in \cite[Ch.~8.7]{q}, we perform a change of variables $z' \to z'-(m+1/2)\delta$ and write $\alpha_{m,m}$ as:
\begin{align}
&\alpha_{m,m} = \int\limits_{z'=-\delta/2}^{\delta/2} dz' r_S(z'+(m+1/2)\delta) \sqrt{1+r_S'^2(z'+(m+1/2)\delta)} \int \limits_{\theta = -\pi}^{\pi} \cos\theta \frac{\exp(-jk_0 R_m)}{4\pi R_m}\ d \theta \nonumber \\
&= 4\int\limits_{z'=0}^{\delta/2} dz' r_S(z'+(m+1/2)\delta) \sqrt{1+r_S'^2(z'+(m+1/2)\delta)}  \int \limits_{\theta = 0}^{\pi} \cos\theta \frac{\exp(-jk_0 R_m)}{4\pi R_m}\ d \theta \nonumber
\end{align}
\begin{align}
&\approx 4 r_{m+1/2}\sqrt{1+s^2_{m}} \int\limits_{z'=0}^{\delta/2} dz' \int \limits_{\theta = 0}^{\pi} \cos\theta \frac{\exp(-jk_0 R_m)}{4\pi R_m}\ d \theta \nonumber \\
& = 4r_{m+1/2}\sqrt{1+s^2_{m}}\bigg[ \int\limits_{z'=0}^{\delta/2} \int \limits_{\theta = 0}^{\pi} \frac{\exp(-jk_0R_m)-1}{R_m}  dz' \ d\theta +\int\limits_{z' = 0}^{\delta/2} \bigg \{ \int \limits_{\theta = 0}^{\pi} \frac{1}{R_m} d\theta \nonumber  \\ & -\frac{1}{r_{m+1/2}} \ln \bigg(\frac{z'}{8r_{m+1/2}}\bigg)\bigg\} \ dz' +\int\limits_{z' = 0}^{\delta/2} \frac{1}{r_{m+1/2}} \ln \bigg(\frac{z'}{8r_{m+1/2}}\bigg) \ dz'\bigg]
\end{align}
The first integral is not singular and can be evaluated numerically. The second integral is also not singular can be rewritten as:
\begin{align}\label{sec_int}
&\int\limits_{z' = 0}^{\delta/2} \bigg \{ \int \limits_{\theta = 0}^{\pi} \frac{1}{R_m} d\theta -\frac{1}{r_{m+1/2}} \ln \bigg(\frac{z'}{8r_{m+1/2}}\bigg)\bigg\} \ dz'  \nonumber \\ &= \int\limits_0^{\delta/2} \frac{2 \ dz'}{\sqrt{4r_{m+1/2}^2+z^{'2}}} \bigg\{ E\bigg( \sqrt{\frac{4r_{m+1/2}^2}{4r_{m+1/2}^2+z^{'2}}}\bigg) + \frac{1}{r_{m+1/2}} \ln \bigg(\frac{z'}{8r_{m+1/2}} \bigg)\bigg\}
\end{align}
where $E(\chi)$ is the complete elliptic integral of the first kind, defined by:
\begin{equation}
E(\chi) = \int\limits_{0}^{\pi/2} \frac{d\beta}{\sqrt{1-\chi \sin^2\beta}}
\end{equation}

\begin{figure}[htpb]
\centering
\includegraphics[scale = 0.3]{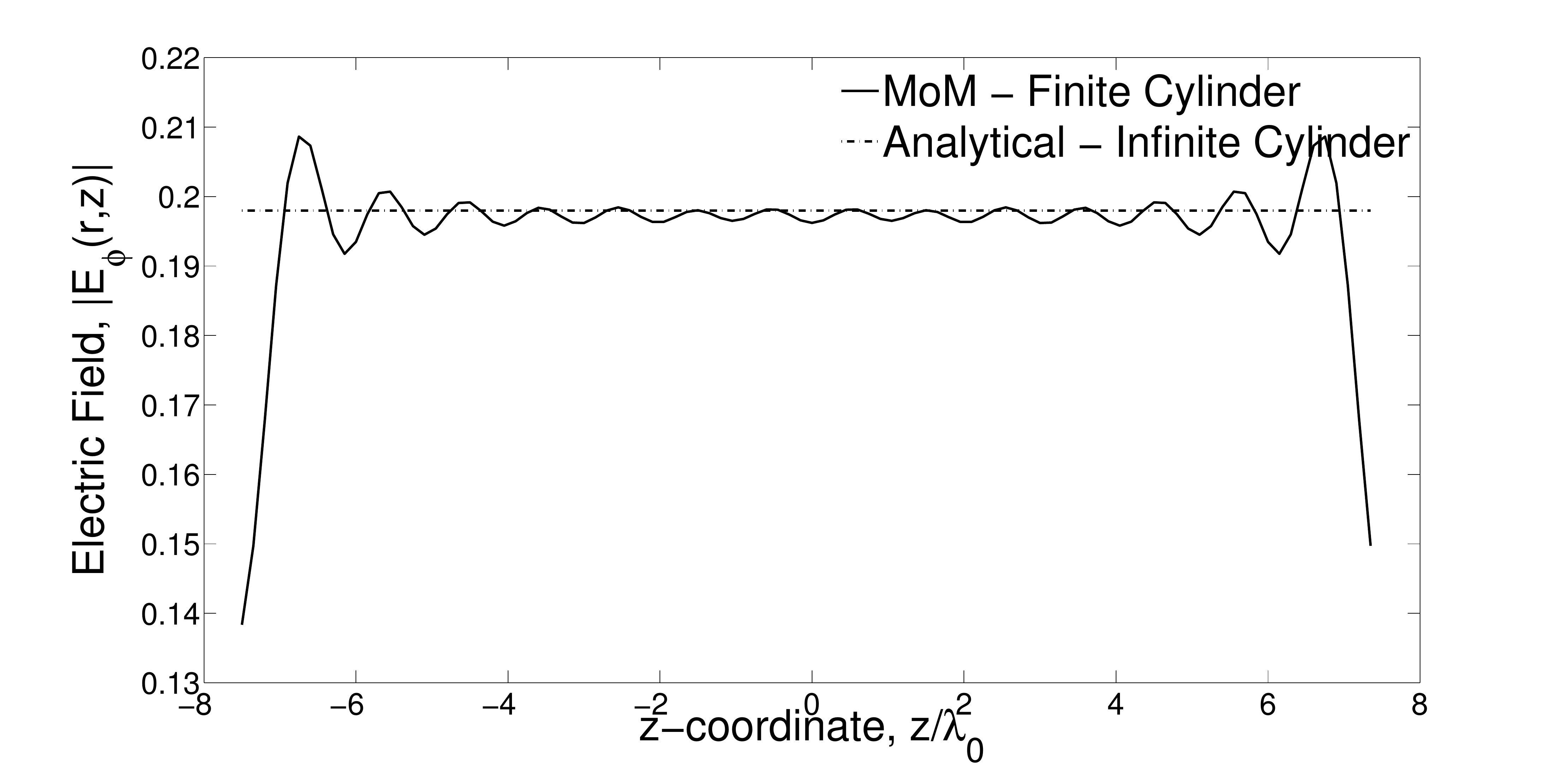}
\caption{Plot of $|E_\phi^{sca}(r,z)|$ as a function of $z$ computed using MoM for an incident field of the form $E_\phi^{inc} = H_1^{(1)}(k_0r)$. In this simulation $a =2\lambda_0$, $r = 2.5\lambda_0$, $\delta = 0.1\lambda_0$ and $L =  15\lambda_0$. The analytically computed fields for an infinitely long cylinder for the same dimensions is also shown to illustrate the diffraction effects due to the finite length of the cylinder.}
\label{Diffraction}
\end{figure}

The integral in Eq.~(\ref{sec_int}) can be evaluated numerically. The last integral can be analytically evaluated to yield:
\begin{equation}
\int\limits_{z' = 0}^{\delta/2} \frac{1}{r_{m+1/2}} \ln \bigg(\frac{z'}{8r_{m+1/2}}\bigg) \ dz' =  \frac{\delta}{2r_{m+1/2}} \bigg[\ln \bigg(\frac{\delta}{16r_{m+1/2}}\bigg)-1  \bigg]
\end{equation}
\end{enumerate}
Finally imposing Eq.~(\ref{red_EFIE}), we obtain:
\begin{equation}
e_m =j\eta_0 k_0 \sum_{n=0}^{N-1} \alpha_{m,n} j_n
\end{equation}
which can be solved to obtain $j_n$. Once $j_n$ are known, the scattered fields at any point $(r,z)$ can be evaluated. For instance, $E_\phi^{sca}(r,z)$ is given by:
\begin{equation}
E_\phi^{sca}(r,z) = -jk_0 \eta_0 A(r,z)
\end{equation}
where
\begin{equation}
A(r,z) = \int\limits_{z'=-\infty}^{\infty} dz' J(z') r_S(z') \sqrt{1+r_S'^2(z')} \int \limits_{\theta = -\pi}^{\pi} \cos\theta \frac{\exp(-jk_0 R)}{4\pi R}\ d \theta
\end{equation}
Substituting for $J(z')$ from Eq.~(\ref{J_appx}), we obtain:
\begin{equation}
E_\phi^{sca}(r,z) =-jk_0\eta_0 \sum_{n=0}^{N-1} \epsilon_n j_n
\end{equation}
where
\begin{align}
\epsilon_n = \int\limits_{z'=n\delta}^{(n+1)\delta} dz' r_S(z') \sqrt{1+r_S'^2(z')} \int \limits_{\theta = -\pi}^{\pi} \cos\theta \frac{\exp(-jk_0 R)}{4\pi R}\ d \theta \nonumber \\
\approx \delta r_{n+1/2} \sqrt{1+s_n^2} \int \limits_{\theta = -\pi}^{\pi} \cos\theta \frac{\exp(-jk_0 R_n)}{4\pi R_n}\ d \theta
\end{align}
with $R_n =\sqrt{r_{n+1/2}^2+r^2-2r r_{n+1/2}\cos\theta +(z-(n+1/2)\delta)^2}$.

 \begin{figure}[htpb]
 \centering
 \begin{subfigure}{\textwidth}
  \centering
\includegraphics[scale = 0.28]{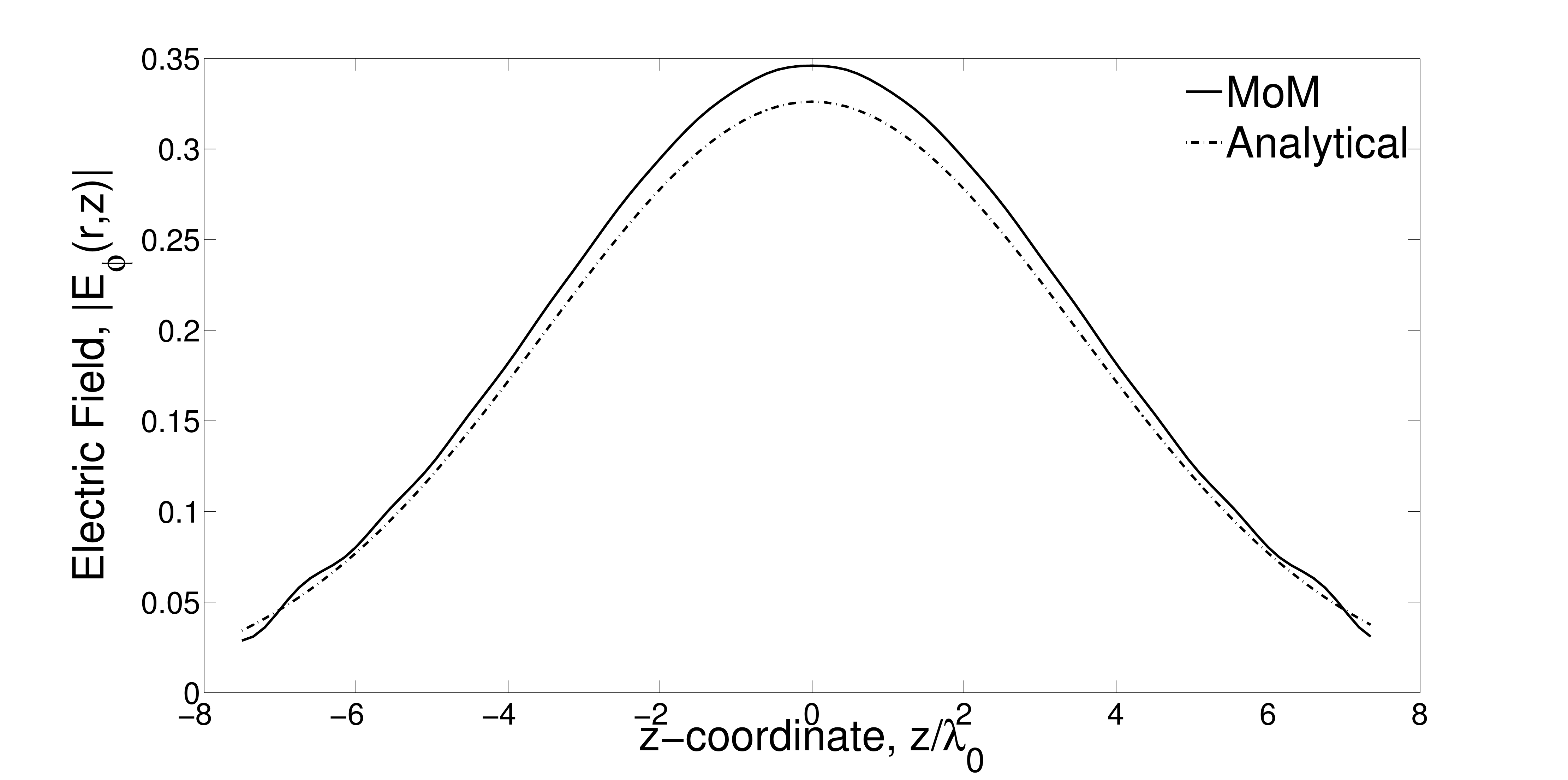}
\caption{$N = 50, \delta = 0.3 \lambda_0$}
\end{subfigure}
\begin{subfigure}{\textwidth}
 \centering
\includegraphics[scale = 0.28]{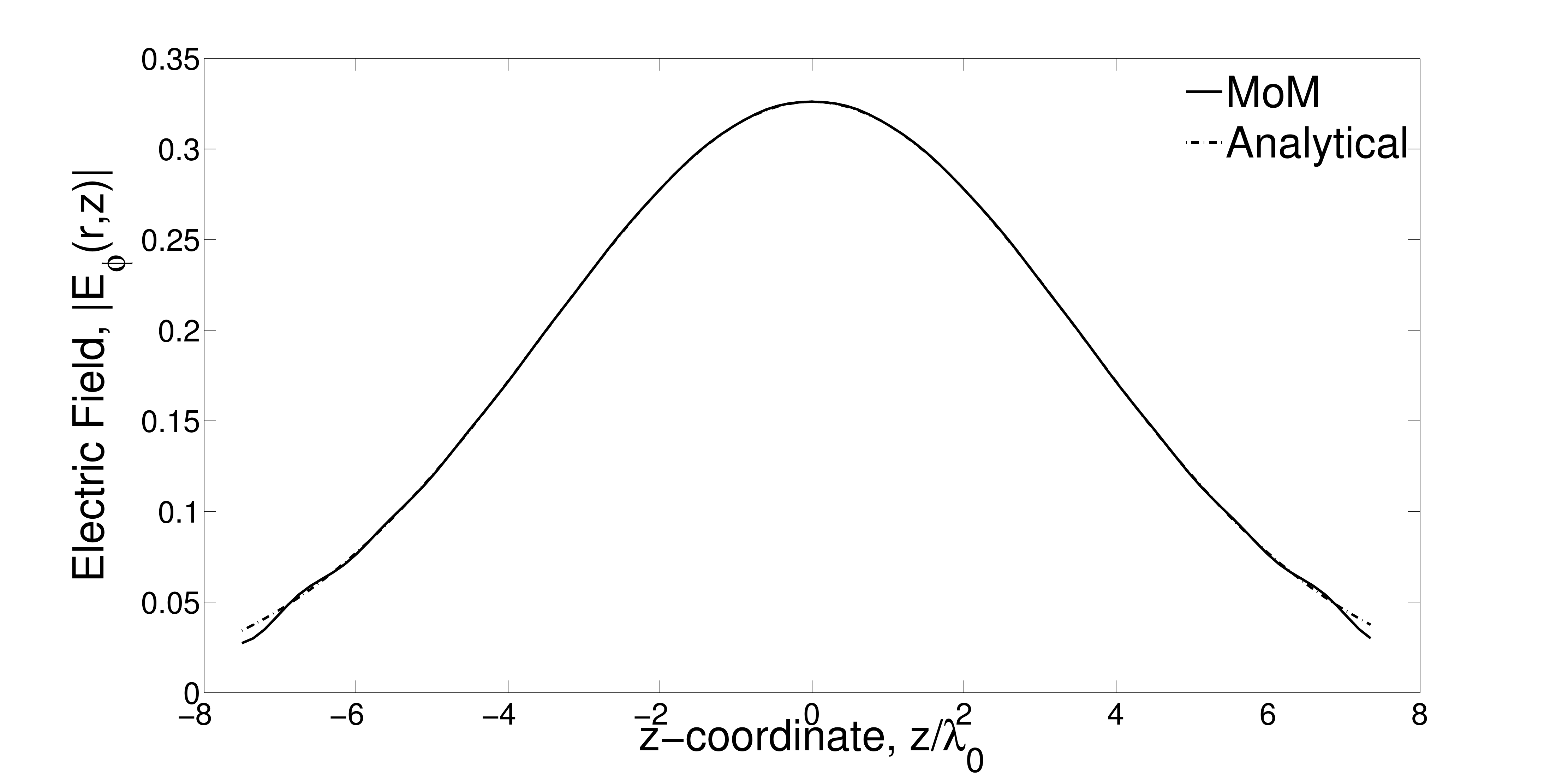}
\caption{$N = 150, \delta = 0.1 \lambda_0$}
\end{subfigure}
\caption{Plot of $|E_\phi^{sca}(r,z)|$ as a function of $z$ computed using MoM and compared with the analytical solution for different discretizations of the cylinder surface. In all the calculations, $a =2\lambda_0$,$r = 2.5\lambda_0$, $w_0 = 5\lambda_0$ and $L =  15\lambda_0$.}
\label{Axvarnh0}
\end{figure}

It is not possible to simulate the response of an infinitely long cylinder with axial roughness, since the MoM would require us to truncate the cylinder in the $z$ direction. Truncation of the cylinder results in diffraction effects from the cylinder edges (as is shown in Fig.~\ref{Diffraction}), which is undesirable if the objective of this exercise is to simulate the response of an \emph{infinitely} long cylinder. However, a standard approach to this problem is to use a tapered incident field \cite{t}. For instance, consider an incident field of the form:
\begin{equation}\label{gauss_inc}
E_{\phi}^{inc} =w_0\int \limits_{-\infty}^{\infty} \exp\bigg(-\frac{k^2 w_0^2}{4}\bigg) H_1^{(1)}(r\sqrt{k_0^2-k^2} ) \exp(jkz) \ dk
\end{equation}
Note that $H_1^{(1)}(r\sqrt{k_0^2-k^2}) \exp(jkz)$ is a cylindrical wave propagating from infinity to the axis of the cylinder and with a component $k\hat{z}$ of the wave-vector along the $z$ direction. By adding such waves with a gaussian envelope $\exp(-k^2 w_0^2/4)$, we have constructed an incident field which is approximately gaussian in the $z$ coordinate. This can be explicitly seen if $k_0 w_0 >> 1$, in which case only very small $k$ contribute significantly to the integral in Eq. (\ref{gauss_inc}). This allows us to make the following Taylor series expansion to the second order in $k$:
\begin{equation}
H_1^{(1)}(r\sqrt{k_0^2-k^2}) \approx H_1^{(1)}(k_0 r) -\frac{k^2 r}{2k_0} H_1^{(1)'}(k_0 r)
\end{equation}
The integral in Eq. (\ref{gauss_inc}) can then be evaluated analytically to yield:
\begin{equation}
E_{\phi}^{inc} (r,z) =2\sqrt{\pi} \exp\bigg(-\frac{z^2}{w_0^2}\bigg) \bigg[H_1^{(1)}(k_0r)+\frac{2rz^2}{k_0w_0^4}H_1^{(1)'}(k_0 r)-\frac{r}{k_0w_0^2} H_1^{(1)'}(k_0r) \bigg]
\end{equation}

Thus, if $L \gg w_0$, then the cylinder behave almost as if it is an infinitely long cylinder as far as the gaussian incident field is concerned. Fig.~\ref{Axvarnh0} shows the comparison between analytical results and MoM results for a smooth cylinder for $N = 50$ and $N = 150$. Clearly, $\delta \sim \frac{\lambda_0}{10}$ is sufficient to accurately simulate the response of the cylindrical scatterer.

\bibliographystyle{IEEEtran}
\bibliography{IEEEabrv,momref}

\end{document}